\documentclass[12pt,fleqn,amssymb,epsfig,a4]{article}
\usepackage{enumerate}
\usepackage{amsmath}
\usepackage{amsfonts}
\usepackage{amssymb}
\usepackage{graphicx}

\newtheorem{theorem}{Theorem}[section]

\newtheorem{lemma}[theorem]{Lemma}
\newtheorem{definition}[theorem]{Definition}

\newtheorem{claim}{Claim}[section]
\newtheorem{conjecture}[theorem]{Conjecture}

\begin{document}

\title{Nowhere-zero 5-flows on cubic graphs with oddness 4}

\author{Giuseppe Mazzuoccolo\thanks{Universit\`{a} di Modena e Reggio Emilia, Dipartimento di Scienze Fisiche, Informatiche e Matematiche, Via Campi 213/b, 41125 Modena, Italy; giuseppe.mazzuoccolo@unimore.it}, Eckhard Steffen\thanks{Paderborn Institute for Advanced Studies in Computer Science and Engineering, Universit\"at Paderborn, 
Warburger Stra{\ss}e 100, 33098 Paderborn, Germany; es@upb.de}}

\date{}

\maketitle

\begin{abstract}
Tutte's 5-Flow Conjecture from 1954 states that every bridgeless graph has a nowhere-zero 5-flow.
In 2004, Kochol proved that the conjecture is equivalent to its restriction on cyclically 6-edge connected cubic graphs.
We prove that every cyclically 6-edge-connected cubic graph with oddness at most 4 has a nowhere-zero 5-flow.

\end{abstract}  

\section[]{Introduction}

An integer nowhere-zero $k$-flow on a graph $G$ is an assignment of a direction and a value of $\{1, \dots, (k-1)\}$ to each edge of $G$ such that the Kirchhoff's law is satisfied at every vertex of $G$. This is the most restrictive definition of a nowhere-zero $k$-flow. But it is equivalent to more flexible definitions, see e.g.~\cite{Seymour_95}.
A cubic graph $G$ is bipartite if and only if it has a nowhere-zero 3-flow, and $\chi'(G)=3$  if and only if $G$ has a nowhere-zero 4-flow. Seymour \cite{Seymour_81} proved that every bridgeless graph has a nowhere-zero 6-flow. 
So far this is the best approximation to Tutte's famous 5-flow conjecture, which is equivalent to its restriction to cubic graphs. 

\begin{conjecture} [\cite{Tutte_54}] \label{5FC}
Every bridgeless graph has a nowhere-zero 5-flow.
\end{conjecture}

Kochol \cite{Kochol_04} proved that a minimum counterexample to the 5-flow conjecture is a cyclically 6-edge-connected cubic graph. Hence, it suffices to prove Conjecture \ref{5FC} for these graphs. 

A classical parameter to measure how far a cubic graph is from being 3-edge-colorable is its {\it oddness}. The oddness, denoted by $\omega(G)$, of a bridgeless cubic graph $G$ is the minimum number of odd circuits in a 2-factor of $G$. Since $G$ has an even number of vertices, $\omega(G)$ is necessarily even. Furthermore, $\omega(G)=0$ if and only if $G$ is $3$-edge-colorable. 

Jaeger \cite{Jaeger_88} showed that cubic graphs with oddness at most 2 have a nowhere-zero 5-flow.
Furthermore, a consequence of the main result in \cite{Steffen_10} is that cyclically $7$-edge connected cubic graphs with oddness at most $4$ have a nowhere-zero 5-flow. 
The following is our main theorem, and it is a strengthening for both previous results.
 
\begin{theorem} \label{oddness4_5flow}
Let $G$ be a cyclically 6-edge-connected cubic graph. If $\omega(G) \leq 4$, then $G$ has a nowhere-zero 5-flow.
\end{theorem}

\section[] {Balanced valuations and flow partitions} \label{gamma_2}

In this section, we recall the concept of flow partitions, which was introduced by the second author in \cite{Steffen_10}.

Let $G$ be a graph and $S \subseteq V(G)$. The set of edges with precisely one end in $S$ is
denoted by $\partial_G(S)$. 
 
An {\em orientation} $D$ of $G$ is an assignment of a
direction to each edge. For $S \subseteq V(G)$, 
$D^-(S)$ ($D^+(S)$) is the set of edges of $\partial_G(S)$ whose head 
(tail) is incident to a vertex of $S$.
The oriented graph is denoted by $D(G)$, 
$d_{D(G)}^-(v) = |D^-(\{v\})|$ and $d_{D(G)}^+(v) = |D^+(\{v\})|$ denote the {\em indegree}
and {\em outdegree} of vertex $v$ in $D(G)$, respectively. The 
degree of a vertex $v$ in the undirected graph $G$ is $d_{D(G)}^+(v) + d_{D(G)}^-(v)$, and it is denoted by $d_G(v)$.  

Let $k$ be a positive integer, and $\varphi$  a function from the edge set of the directed graph $D(G)$
into the set $\{0, 1, \dots, k-1\}$. For $S \subseteq V(G)$ let 
$\delta \varphi (S) = \sum_{e \in D^+(S)}\varphi(e) - \sum_{e \in D^-(S)}\varphi(e)$. 
The function $\varphi$ is a $k$-flow on $G$ if $\delta \varphi(S) = 0$ for every $S \subseteq V(G)$. 
 The {\em support} of $\varphi$  is 
the set $\{e \in E(G) : \varphi(e) \not = 0\}$, and it is denoted by $supp(\varphi)$.
A $k$-flow $\varphi$ is a nowhere-zero $k$-flow if $supp(\varphi) = E(G)$.

We will use balanced valuations of graphs, which were introduced by 
Bondy \cite{Bondy} and Jaeger \cite{Jaeger_75}. A {\em balanced valuation} of a graph
$G$ is a function $f$ from the vertex set $V(G)$ into the real numbers, such that
$| \sum_{v \in X} f(v) | \leq | \partial_G(X) |$ for all $X \subseteq V(G)$. Jaeger proved 
the following fundamental theorem. 

\begin{theorem} [\cite{Jaeger_75}] \label{Thm_Jaeger_75} 
Let $G$ be a graph with orientation $D$ and $k\geq 3$. Then $G$ 
has a nowhere-zero $k$-flow if and only if there 
is a balanced valuation $f$ of $G$ with 
$ f(v) = \frac{k}{k-2}(2d_{D(G)}^+(v) - d_G(v))$, for all $v \in V(G).$
\end{theorem}

 In particular, Theorem \ref{Thm_Jaeger_75} says that a cubic graph $G$ has a nowhere-zero $5$-flow 
 if and only if there is a balanced valuation of $G$ with values in $\{ \pm \frac{5}{3}\}$.

Let $G$ be a bridgeless cubic graph, and ${\cal F}_2$ be a $2$-factor of $G$ with odd circuits $C_1,\ldots,C_{2t}$, 
and even circuits $C_{2t+1},\ldots, C_{2t+l}$ ($t \geq 0$, $l \geq 0$), and let ${\cal F}_1$ be the complementary $1$-factor.

A {\em canonical} $4$-edge-coloring, denoted by $c$, of $G$ with respect to ${\cal F}_2$ colors the edges of ${\cal F}_1$ with color $1$, the edges of the even circuits of ${\cal F}_2$ with $2$ and $3$, alternately, and the edges of the odd circuits of ${\cal F}_2$ with colors $2$ and $3$ alternately, but one edge which is colored $0$. Then, there are precisely $2t$ vertices $z_1,\ldots,z_{2t}$ where color $2$ is missing (that is, no edge which is incident to $z_i$ has color $2$).

The subgraph which is induced by the edges of colors $1$ and $2$ is union of even circuits and $t$ paths $P_i$ of odd length 
and with $z_1,\ldots,z_{2t}$ as ends.
Without loss of generality we can assume that $P_i$ has ends $z_{2i-1}$ and $z_{2i}$, for $i \in \{1,\ldots,t\}$.

Let $M_G$ be the graph obtained from $G$ by adding two edges $f_i$ and $f'_i$ between $z_{2i-1}$ and $z_{2i}$ for $i \in \{1,\ldots,t\}$. 
Extend the previous edge-coloring to a proper edge-coloring of $M_G$ by coloring $f'_i$ with color $2$ and $f_i$ with color $4$.
Let $C'_1,\ldots,C'_s$ be the cycles of the $2$-factor of $M_G$ induced by the edges of colors $1$ and $2$ ($s \geq t$). 
In particular, $C'_i$ is the even circuit obtained by adding the edge $f'_i$ to the path $P_i$, for $i \in \{1,\ldots,t\}$.
Finally, for $i \in \{1,\ldots,t \}$ let $C''_i$ be the $2$-circuit induced by the edges $f_i$ and $f'_i$. 
We construct a nowhere-zero $4$-flow on $M_G$ as follows:
\begin{itemize}
\item for $i \in \{1, \dots, 2t+l\}$ let $(D_i,\varphi_i)$ be a nowhere-zero flow on the directed circuit $C_i$ with $\varphi_i(e)=2$ for all $e \in E(C_i)$;
\item for $i \in \{1 \dots, s\}$ let $(D'_i,\varphi'_i)$ be a nowhere-zero flow on the directed circuit $C'_i$ with $\varphi'_i(e)=1$ for all $e \in E(C'_i)$; 
\item for $i \in \{1, \dots, t\}$ let $(D''_i,\varphi''_i)$ be a nowhere-zero flow on the directed circuit $C''_i$ (choose $D''_i$ such that $f'_i$ receives the same direction as in $D'_i$) with $\varphi''_i(e)=1$ for all $e \in \{f_i,f'_i\}$.
\end{itemize}
Then,

$$(D,\varphi)= \sum_{i=1}^{2t+l} (D_i,\varphi_i) + \sum_{i=1}^{s} (D'_i,\varphi'_i) + \sum_{i=1}^{t} (D''_i,\varphi''_i) $$

is the desired nowhere-zero $4$-flow on $M_G$. 

By Theorem \ref{Thm_Jaeger_75}, there is a balanced valuation $w(v)=2(2d_{D(M_G)}^+(v) - d_{M_G}(v))$ of $M_G$. It holds that $|2d_{D(M_G)}^+(v) - d_{M_G}(v)|=1$, and hence, $w(v) \in \{ \pm 2 \}$ for all vertices $v$. 
The vertices of $M_G$, and therefore, of $G$ as well, are partitioned into two classes $A=\{v | w(v)=-2\}$ and $B=\{v | w(v)=2\}$. 
We call the elements of $A$ ($B$) the white (black) vertices of $G$, respectively.

\begin{definition}
Let $G$ be a bridgeless cubic graph and ${\cal F}_2$ a 2-factor of $G$. A partition of $V(G)$ into two classes $A$ and $B$ constructed as above with a canonical $4$-edge-coloring $c$, the $4$-flow $(D,\varphi)$ on $M_G$ and the induced balanced valuation $w$ of $M_G$ is called a \textbf{flow partition} of $G$ w.r.t.~${\cal F}_2$. The partition is 
denoted by $P_G(A,B) (=P_G(A,B,{\cal F}_2,c,(D,\varphi),w))$.  
\end{definition}

\begin{lemma}  \label{diffends_lemma} Let $G$ be a bridgeless cubic graph and
$P_G(A,B)$ be a flow partition of $V(G)$ which is induced by a canonical nowhere-zero 4-flow with respect to an edge-coloring $c$. Let $x$, $y$ be the two vertices of an edge $e$. If $e \in c^{-1}(1) \cup c^{-1}(2)$, then $x$ and $y$ belong to different classes, i.e. $x \in A$ if and only if $y \in B$.
\end{lemma}

From a flow partition $P_G(A,B) (=P_G(A,B,{\cal F}_2,c,(D,\varphi),w))$ we easily obtain a flow partition 
$P_G(A',B') (=P_G(A',B',{\cal F}_2,c,(D',\varphi'),w'))$ such that the colors on the vertices of $P_i$ are switched. 
Let $(D',\varphi')$ be the nowhere-zero $4$-flow on $M_G$ obtained by using the same $2$-factor ${\cal F}_2$, the same $4$-edge-coloring $c$ of $G$ and the same orientations for all circuits, but for one $i \in \{i, \dots ,t\}$ use opposite orientation 
of $C_i'$ and $C_i''$ with respect to the one selected in $(D,\varphi)$. 

\begin{lemma}\label{switching_lemma} Let $G$ be a bridgeless cubic graph and 
$P_G(A,B)$ be the flow partition which is  induced by the nowhere-zero $4$-flow $(D,\varphi)$. If $P_G(A',B')$ is the flow partition induced by the nowhere-zero $4$-flow $(D',\varphi')$, then 
 $A \setminus V(P_i) = A' \setminus V(P_i)$,  
$B \setminus V(P_i) = B' \setminus V(P_i)$,  $A \cap V(P_i) = B' \cap V(P_i)$ and $B \cap V(P_i) = A' \cap V(P_i)$.
\end{lemma}

\section{Proof of Theorem \ref{oddness4_5flow}}

Suppose to the contrary that the statement is not true. Then there is a cyclically 6-edge-connected cubic graph $G$, which 
has no nowhere-zero 5-flow. Let ${\cal F}_2$ be a 2-factor of $G$ with
precisely four odd circuits $C_1,\dots,C_4$. Let $c$ be a canonical 4-edge coloring of $G$ and 
$z_1, z_2, z_3, z_4$  be the four vertices where color $2$ is missing. Let $Z=\{z_1$,$z_2$,$z_3,z_4\}$.
Note, that in any flow partition which depends on ${\cal F}_2$ and $c$, the vertices 
$z_1$ and $z_2$ (and $z_3$ and $z_4$ as well) belong to different color classes. 
By Lemma \ref{switching_lemma} there are flow partitions $P_G(A,B)$ and $P_G(A',B')$ of $G$ 
such that $\{z_1, z_3\} \subseteq A$, and $\{z_1, z_4\} \subseteq A'$. Hence, $\{z_2, z_4\} \subseteq B$ 
and $\{z_2, z_3\} \subseteq B'$.

Let $w$ be the function with $w(v)=-\frac{5}{3}$ if $v \in A$ and $w(v)=\frac{5}{3}$ if $v \in B$, and 
 $w'$ be a function with $w'(v)=-\frac{5}{3}$ if $v \in A'$ and $w'(v)=\frac{5}{3}$ if $v \in B'$.

We will prove that $w$ or $w'$ is a balanced valuation of $G$, and therefore, $G$ has a nowhere-zero $5$-flow by Theorem \ref{Thm_Jaeger_75}. Hence, there is no counterexample and Theorem \ref{oddness4_5flow} is proved.

\subsection{$Z$-separating edge-cuts}

Since $G$ has no nowhere-zero 5-flow, $w$ and $w'$ are not balanced valuations of $G$. 
Then there are $S\subseteq V(G)$, $S' \subseteq V(G)$ with 
$|\sum_{v \in S} w(v)| > | \partial_G(S) |$, and $|\sum_{v \in S'} w'(v)| > | \partial_G(S') |.$

We will prove some properties of the edge-cuts $\partial_G(S)$ and $\partial_G(S')$. We deduce the results
for $S$ only. The results for $S'$ follow analogously. 
If $S=V(G)$, then $|\sum_{v \in S} w(v)|= 0 = | \partial_G(S) |$. 
Therefore, $S$, $S'$ are a proper subset of $V$.
If $S = \{v\}$, then  $|\sum_{v \in S} w(v)|= \frac{5}{3} \leq 3 = |\partial_G(S)|$. 
Since $G$ is cyclically $6$-edge-connected, it has no non-trivial $3$-edge-cut and no $2$-edge-cut. Hence, 
we assume that $|\partial_G(S)| \geq 4$ in the following. 

Let $k$ ($k'$)  be the absolute value of the difference between the number of black and white vertices in $S$ ($S'$). Hence, 
$\frac{5}{3} k > |\partial_G(S)|$, and $\frac{5}{3} k' > |\partial_G(S')|$.

For $i \in \{0,1,2,3\}$, let $c_i=|\partial_G(S) \cap c^{-1}(i)|$ and $c'_i=|\partial_G(S') \cap c^{-1}(i)|$.

\begin{claim} \label{c1}$|\partial_G(S)| \equiv k \pmod 2$,  $|\partial_G(S')| \equiv k' \pmod 2$
\end{claim}
{\it Proof.} If $k$ is even, then $|S \cap A|$ and $ |S \cap B|$ have the same parity, and if $k$ is
odd, then they have different parities. Since $S$ is the disjoint union of $S \cap A$ 
and $S \cap B$ it follows that $k$ and $|S|$ have the same parity. Since $G$ is cubic it
follows that $|\partial_G(S)| \equiv k \pmod 2$. $\square$

Let $q_A$ ($q_B$) be the number of white (black) vertices of $S$ where color $2$ is missing. Let $q=|q_A - q_B|$. 
Since $Z$ has two black and two white vertices, it follows that $q \leq 2$. 

\begin{claim} \label{c2}$|S \cap Z|= 2 = q$, and $|S' \cap Z|=2=q'$.
\end{claim}
{\it Proof.}  
Since $c^{-1}(1)$ is a $1$-factor of $G$, Lemma \ref{diffends_lemma} implies that $k = c_1$. Hence, 
\begin{equation}\label{eqn_c1}
c_1 > \frac{3}{5} |\partial_G(S)|.
\end{equation}

Furthermore, Lemma \ref{diffends_lemma} implies that $k \leq c_2 + q$. Hence, 
\begin{equation}\label{eqn_c2}
c_2+q > \frac{3}{5} |\partial_G(S)|. 
\end{equation}

Suppose to the contrary, that $|S \cap Z|\not=2$. Thus, $q \leq 1$ and $c_2+1 \geq c_1$.
Hence, $|\partial_G(S)| \geq c_1 + c_2 \geq 2k-1$, and therefore, $\frac{5}{3}k \leq |\partial_G(S)|$ if 
$|\partial_G(S)| \geq 6$. 

If $|\partial_G(S)|=4$, then $k \leq 2$, and  if $|\partial_G(S)|=5$, then $k \leq 3$.
In both cases, it follows that $\frac{5}{3}k \leq |\partial_G(S)|$, a contradiction. 
Thus, $|S \cap Z|=2$, and therefore, $q \in \{0,2\}$. If $q = 0$, then $|\partial_G(S)| \geq c_1+c_2 \geq 2k$, a contradiction.
Hence, $q = 2$. $\square$ 

\begin{claim} \label{c3} $|\partial_G(S)|=6$, $c_1=4$ and $c_2=2$, and $|\partial_G(S')|=6$, $c'_1=4$ and $c'_2=2$.
\end{claim}
{\it Proof.}  
If $|\partial_G(S)| = 4$, then $k \geq 3$. Hence, $c_1 = 3$ and $c_2 =1$. The edge of $\partial_G(S) \cap c^{-1}(2)$
is contained in a circuit of ${\cal F}_2$ whose edges are not in $c^{-1}(1)$. Hence, $2 \geq c_1 = k$, a contradiction. 
If $|\partial_G(S)| = 5$, then $k \geq 4$. We deduce a contradiction as in the case before.  

Now suppose to the contrary that $|\partial_G(S)|> 6$.
Since  $c_1 > \frac{3}{5} |\partial_G(S)|$, $c_2 >  \frac{3}{5} |\partial_G(S)| - 2$, and $c_1 + c_2 \leq |\partial_G(S)|$, it follows
that $|\partial_G(S)| > \frac{6}{5} |\partial_G(S)| - 2 $. Therefore, $|\partial_G(S)|<10$.
If $|\partial_G(S)|=7$, then $c_1 \geq 5$ and $c_2 \geq 3$, a contradiction.
If $|\partial_G(S)|=8$, then $c_1 = 5$ and $c_2 = 3$, a contradiction to Claim \ref{c1} since $c_1 = k$.
If $|\partial_G(S)|=9$, then $c_1 \geq 6$ and $c_2 \geq 4$, a contradiction.
Hence, $|\partial_G(S)|=6$ and $c_1 \geq 4$ and $c_2 \geq 2$. That leaves the unique possibility $c_1=4$ and $c_2=2$. $\square$

\begin{claim} \label{c4}
$G[S]$ and $G[S']$ are connected.
\end{claim}
{\it Proof.} If $G[S]$ is not connected, then there is $E \subseteq \partial_G(S)$ such that
$G-E$ has at least two components $K_1$ and $K_2$. Since $G$ does not have a 2-edge-cut or a non-trivial 3-edge-cut, 
it follows that $|E|=3$ and one of $K_1$, $K_2$ is a single vertex. Hence, $\frac{5}{3}k \leq |\partial_G(S)|$, a contradiction. $\square$

\begin{definition} \label{D_bad} A 6-edge-cut $E$ of $G$ is \textbf{bad} with respect to a flow partition $P_G(A^*,B^*)$
if it satisfies the following two conditions: 
\begin{enumerate}[i)]
 \item $|E \cap c^{-1}(1)| = 4$ and $|E \cap c^{-1}(2)| = 2$,
 \item $E$ partitions the vertices $z_1$, $z_2$, $z_3$ and $z_4$ into two sets $\{z_{i_1},z_{i_2}\}$, $\{z_{i_3},z_{i_4}\}$, which are in 
different components of $G-E$ and $\{z_{i_1},z_{i_2}\} \subseteq A^*$ or $\{z_{i_1},z_{i_2}\} \subseteq B^*$.
\end{enumerate}    
\end{definition}

Note that $\{z_{i_1},z_{i_2}\} \subseteq A^*$ if and only if $\{z_{i_3},z_{i_4}\} \subseteq B^*$. Further,
only condition $ii)$ depends on the flow partition. Condition i) depends on the canonical 4-edge-coloring of 
$G$ which is unchanged along the proof. From the previous results we deduce:

\begin{claim} \label{c5}
$\partial_G(S)$ is bad w.r.t.~$P_G(A,B)$ and $\partial_G(S')$ is bad w.r.t.~$P_G(A',B')$. 
\end{claim}

Bad 6-edges-cuts are the only obstacles in $G$ for having a nowhere-zero 5-flow. 
In order to deduce the desired contradiction we will show that 
all 6-edge-cuts are not bad with respect to either $P_G(A,B)$ or $P_G(A',B')$.

Recall that, $z_1$ and $z_3$  receive the same color in $P_G(A,B)$, and that $z_1$ and $z_4$ receive the same color in $P_G(A',B')$.
For $i \in \{1,2,3\}$, let ${\cal S}_i = \{V : V \subseteq V(G) \mbox{ and } \{z_1,z_i\} \subseteq V\}$ and 
${\cal E}_i = \{E: E \subseteq E(G), V \in {\cal S}_i  \mbox{ and } E = \partial_G(V) \}$ 
be the corresponding set of edge-cuts. 
Since $z_1$ and $z_2$ have different colors in both $P_G(A,B)$ and $P_G(A',B')$,  all edge-cuts in ${\cal E}_2$ 
are not bad with respect to $P_G(A,B)$ and with respect to $P_G(A',B')$.

\begin{figure}
\centering
\includegraphics[width=8cm]{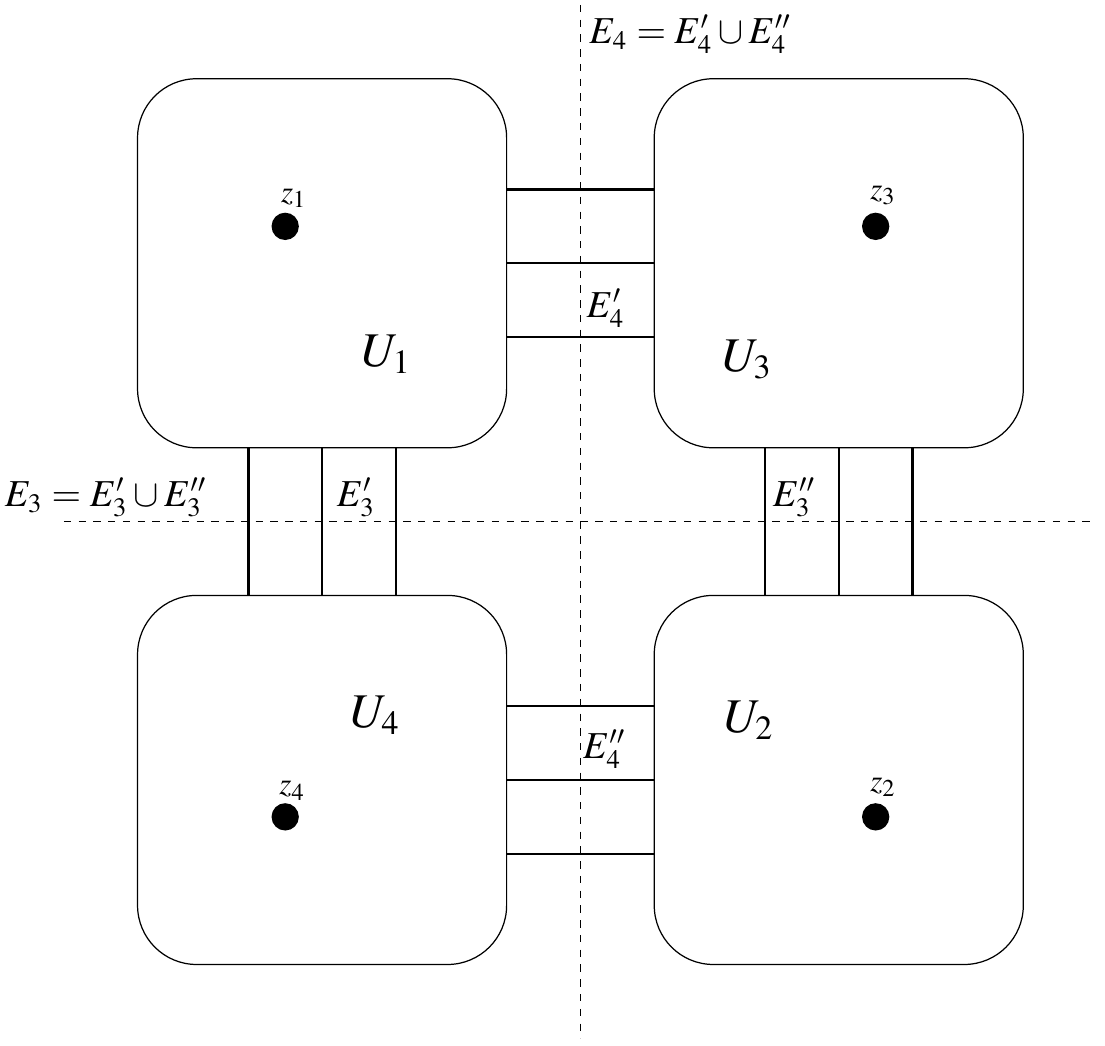}
\caption{$Z$-separating 6-edge cuts}\label{badcuts}
\end{figure}

For $i \in \{3,4\}$, by Claim \ref{c5} there is a 6-edge-cut $E_i \in {\cal E}_i$ which is bad. 
By Claim \ref{c4}, $G-E_3$ consists of two components with vertex sets $X$ and $Y$, i.e.~$X \cup Y=V(G)$.
Analogously, $G-E_4$ consists of  two components with vertex sets $X'$ and $Y'$. 
Let $U_1=X \cap X'$, $U_2=Y \cap Y'$, $U_3=X \cap Y'$ and $U_4=Y \cap X'$. 
Thus, $z_i \in U_i$ for $i \in \{1,\dots,4\}$, see Figure \ref{badcuts}.

\begin{claim}\label{small_component} 
$|\partial_G(U_i)| \geq 5$. In particular, $|\partial_G(U_i)|=5$ if and only if 
$G[U_i]$ is a path with two edges, one of color $0$ and one of color $3$. 
 \end{claim}
{\it Proof.} If $G[U_i]$ has a circuit, then $|\partial_G(U_i)| \geq 6$ since $G$ is cyclically $6$-edge-connected.
If this is not the case, then $G[U_i]$ is a forest, say with $n$ vertices. Hence, $|\partial_G(U_i)| \geq n+2$.
Since $\partial_G(U_i) \subseteq E_3 \cup E_4$, 
it follows that $\partial_G(U_i) \subseteq c^{-1}(1) \cup c^{-1}(2)$.
Two edges $z_ix_i$ and  $z_iy_i$ which are incident to $z_i$ are colored with color 0 and 3, respectively.
Hence, $\{x_i, y_i\} \subseteq U_i$, $n \geq 3$, and $\partial_G(U_i) \geq 5$.
If $|\partial_G(U_i)| = 5$, then $|U_i|= 3$, and $G[U_i]$ is a path with two edges, one of color $0$ and one of color $3$.
$\square$

\begin{claim}\label{two_large_components}
$|\partial_G(U_i)| = 5$ for at most two of the four subsets $U_i$. Furthermore,
if there are $i,j$ such that $i \not = j$ and $|\partial_G(U_i)| = |\partial_G(U_j)| = 5$, 
then  $\{i,j\} \in \{\{1,2\}, \{3,4\}\}$. 
\end{claim}
{\it Proof.} Since $E_3$ and $E_4$ are bad, each of them has exactly two edges of color $2$ and four edges of color $1$. 
Hence, each of them intersects with at most one circuit of ${\cal F}_2$. For each $i \in \{1, \dots ,4\}$, $|U_i \cap c^{-1}(0)| = 1$, and
hence, there are $j_1, j_2$ such that $j_1 \not= j_2$ and $U_{j_1}$, $U_{j_2}$ contain an odd circuit of ${\cal F}_2$. 
Since $G$ is cyclically 6-edge-connected it follows that $|\partial_G(U_{j_1})| \geq 6$ and $|\partial_G(U_{j_2})| \geq 6$.

Let $i,j \in \{1, \dots, 4\}$ such that $i \not = j$ and $|\partial_G(U_i)| = |\partial_G(U_j)| = 5$.
For symmetry, it suffices to prove that $\{i,j\}| \not= \{1,3\}$.
Suppose to the contrary that $\{i,j\} = \{1,3\}$. By Claim \ref{small_component}, $G[U_1]$ and $G[U_3]$ are paths of length two with edges colored $0$ and $3$. Further, $\partial_G(U_1)$ consists of three edges of color $1$ and two edges of color $2$, which belong to the
odd circuit $C_1$ of ${\cal F}_2$. Analogously, the two edges of color 2 of $\partial_G(U_3)$ belong to the odd circuit $C_3$ of ${\cal F}_2$.
Hence, both pairs of edges of color $2$ in $\partial_G(U_1)$ and $\partial_G(U_3)$ belong to $E_3$ and they are distinct, a contradiction since $E_3$ has only two edges of color $2$.
$\square$

For $i \not = j$ let $\partial_G(U_i,U_j)$ be the set of edges with one vertex in $U_i$ and the other one in $U_j$.

\begin{claim}
The following relations hold:
\begin{itemize}
 \item $|\partial_G(U_i,U_j)|=0$, for $\{i,j\} \in \{\{1,2\},\{3,4\}\}$.
 \item $|\partial_G(U_i,U_j)|=3$, for $\{i,j\} \in \{\{1,3\}, \{1,4\}, \{2,3\}, \{2,4\}\}$.
\end{itemize}
\end{claim}
{\it Proof.} Recall that  $|E_3|=|E_4|=6$. Hence, $|E_3 \cup E_4| \leq 12$. Due to Claim \ref{two_large_components}, we can assume 
that $|\partial_G(U_1)| \geq 5$, $|\partial_G(U_2)| \geq 5$, $|\partial_G(U_3)|\geq 6$ and $|\partial_G(U_4)|\geq 6$. By adding up, we obtain $\sum_{i=1}^4 |\partial_G(U_i)| \geq 22$, where each edge of $E_3$ and $E_4$ is counted exactly twice. Hence, $|E_3 \cup E_4| \geq 11$. 
If $|E_3 \cup E_4| = 11$, then exactly one edge, say $e$, belongs to $E_3 \cap E_4$. 
If $e \in \partial_G(U_1,U_2)$, then $\partial_G(U_3)$ and $\partial_G(U_4)$ are distinct sets of cardinality at least $6$. Hence, 
$|E_3 \cup E_4| > 12$, a contradiction. 
If $e \in \partial_G(U_3,U_4)$, then $\partial_G(U_1,U_4)$ or $\partial_G(U_2,U_3)$ has cardinality at most 2, say, without loss of generality, 
$\partial_G(U_1,U_4)$. For the same reason, $\partial_G(U_1,U_3)$ or $\partial_G(U_2,U_4)$ has cardinality at most 2.
If $|\partial_G(U_1,U_3)| \leq 2$, then $|\partial_G(U_1)|\leq 4$, and if  $|\partial_G(U_2,U_4)| \leq 2$, then $|\partial_G(U_4)|\leq 5$,
a contradiction (in both cases). 
Hence, $|E_3 \cup E_4| = 12$, and therefore, $|\partial_G(U_i,U_j)|=0$ for $\{i,j\} \in \{\{1,2\}, \{3,4\}\}$. 

Now,   
$|\partial_G(U_i,U_j)|=3$, for $\{i,j\} \in \{\{1,3\}, \{1,4\}, \{2,3\}, \{2,4\}\}$ can be deduced easily. $\square$

Let $E'_3=E_3 \cap \partial(U_1)$, $E''_3=E_3 \cap \partial(U_2)$, and $E'_4=E_4 \cap \partial(U_1)$, $E''_4=E_4 \cap \partial(U_2)$, see Figure \ref{badcuts}.

Let $H = G[c^{-1}(1) \cup c^{-1}(2)]$.  The components of $H$ are even circuits and the two paths $P_1$ and $P_2$,
where $P_1$ has the end vertices $z_1$, $z_2$, and $P_2$ has the end vertices $z_3$, $z_4$.  
The paths $P_1$ and $P_2$ intersect both $E_3=E'_3 \cup E''_3$ and  $E_4=E'_4 \cup E''_4$ an odd number of times, since both,
$E_3$ and $E_4$, separate their ends.  
For symmetry, we can assume that $P_1 \cap E'_3$ and $P_1 \cap E''_4$ are even, and hence, $P_1 \cap E''_3$ and $P_1 \cap E'_4$ are odd. 
Furthermore, we assume that $P_2 \cap E''_3$ and $P_2 \cap E''_4$ are even, and hence, $P_2 \cap E'_3$ and $P_2 \cap E'_4$ are odd. 
Note, that every other possible choice produces an analogous configuration.
The $6$-edge-cut $E'_3 \cup E'_4$ contains an odd number of edges of $E(P_1) \cup E(P_2)$. Since $E'_3 \cup E'_4 \subseteq E(H)$,
it follows that an odd number of edges of $E'_3 \cup E'_4$ are not in $E(P_1) \cup E(P_2)$, a contradiction, since 
all other components of $H$ are circuits, and they intersect every edge-cut an even number of times. 

Hence, at least one of $E_3$ and $E_4$ is not bad, contradicting our assumption that both of them are bad.

\end{document}